\documentclass[12pt,centertags,oneside]{amsart}
\usepackage{amsmath,amstext,amsthm,amscd,typearea,hyperref,stmaryrd}
\usepackage{amssymb}
\usepackage{a4wide}
\usepackage[mathscr]{eucal}
\usepackage{mathrsfs}
\usepackage{typearea}
\usepackage{charter}
\usepackage{pdfsync}
\usepackage[a4paper,width=16.2cm,top=3cm,bottom=3cm]{geometry}

\numberwithin{equation}{section}

\newtheorem{theorem}{Theorem}[section]

\newtheorem{proposition}[theorem]{Proposition}
\newtheorem{corollary}[theorem]{Corollary}
\newtheorem{lemma}[theorem]{Lemma}
\newtheorem{remark}[theorem]{Remark}

\newcommand{\cali}[1]{\mathscr{#1}}

\newcommand{\scl}{{\rm sc}}

\newcommand{\Prob}{\mathop{\mathrm{Prob}}\nolimits}

\newcommand{\Leb}{\mathop{\mathrm{Leb}}\nolimits}

\newcommand{\dist}{\mathop{\mathrm{dist}}\nolimits}

\newcommand{\loc}{{loc}}
\newcommand{\ddc}{dd^c}

\newcommand{\ddbar}{\partial\overline\partial}

\newcommand{\Sym}{{\rm Sym}}

\newcommand{\sign}{{\rm sign\ \!}}

\newcommand{\Cc}{\cali{C}}

\newcommand{\Ec}{\cali{E}}
\newcommand{\Fc}{\cali{F}}
\newcommand{\Gc}{\cali{G}}

\renewcommand{\Mc}{\cali{M}}

\newcommand{\Pc}{\cali{P}}

\newcommand{\FS}{{\rm FS}}

\newcommand{\C}{\mathbb{C}}
\newcommand{\D}{\mathbb{D}}
\newcommand{\E}{\mathbb{E}}

\newcommand{\R}{\mathbb{R}}
\renewcommand\P{\mathbb{P}}

\renewcommand{\S}{\mathbb{S}}

\title[Large deviation theorem for zeros of polynomials]{Large deviation theorem for zeros of polynomials\\ and Hermitian random  matrices}

\author{Tien-Cuong Dinh}
\address{Department of Mathematics, National University 
of Singapore, 10 Lower Kent Ridge Road, Singapore 119076.}
\email{matdtc@nus.edu.sg}

\begin{document}

\begin{abstract}
We give abstract versions of the large deviation theorem for the distribution of zeros of polynomials and apply them to the characteristic polynomials of Hermitian random matrices. We obtain new estimates related to the local semi-circular law for the empirical spectral distribution of these matrices when the 4th moments of their entries are controlled.
\end{abstract}

\maketitle

\medskip

\noindent {\bf Classification AMS 2010 :} 12B52, 60B20. 

\medskip

\noindent {\bf Keywords :} logarithmic potential, large deviation theorem, Hermitian random matrix, semi-circular law.



\section{Introduction}

Let $\mu_0$ be a given probability measure whose support is contained in a smooth compact non-closed curve $K$ in $\C$, that is, $K$ is parametrised by an $\Cc^\infty$ injective map from $[0,1]$ to $\C$ with non-vanishing derivative. Denote by $u_0$ its logarithmic potential which is defined  for $z\in\C$ by 
$$u_0(z):=\int_{a\in\C} \log|z-a| d\mu_0(a).$$
We assume that the function $u_0$ is H\"older continuous on $K$. This property is satisfied when $\mu_0(I)\leq c|I|^\alpha$ for all arcs $I$ in $K$ of length $|I|$, and for some constants $c>0$ and $\alpha>0$. 
Let $\Pc_0^n$ denote the set of all monic complex polynomials of one variable and degree $n$. Consider, for each $n\geq 1$, a probability measure on $\Pc_0^n$. 

The purpose of this paper is to give a sufficient condition so that the zeros of generic polynomials $Q\in\Pc_0^n$ are equidistributed with respect to the measure $\mu_0$, as $n$ tends to infinity. We will give some estimates which are similar to  the classical large deviation theorem in the theory of probability.
The genericity of polynomials is  with respect to the given probability measures on $\Pc_0^n$. We then apply the result to the characteristic polynomials of Hermitian random matrices and obtain new estimates for the distribution of the eigenvalues of these matrices. In this case, the probability measure $\mu_0$ will be the so-called semi-circular distribution, that we will denote by $\mu_\scl$. 

If $Q$ is a polynomial of degree $n$ and $a_1,\ldots, a_n$ are its zeros, define 
$$\mu_Q:={1\over n} \sum_{k=1}^n \delta_{a_k},$$
where $\delta_{a_k}$ is the Dirac mass at the point $a_k$. This is the probability measure equidistributed on the zero set of $Q$. 
Here is our first main theorem with a quite simple hypothesis.

\begin{theorem} \label{th_abstract}
Let $\mu_0$ be a given probability measure whose support is contained in a smooth compact non-closed curve $K$ in $\C$. Denote by $u_0$ its logarithmic potential and assume that $u_0$ is H\"older continuous on $K$. Let $\alpha>0$ be a constant.
Consider a probability measure on the set $\Pc^n_0$ of all monic polynomials of degree $n$ and a constant $c_n\geq 1$  
such that the expectation of $|Q(z)|^\alpha$   for $Q\in \Pc^n_0$ satisfies
$$\E(|Q(z)|^\alpha)\leq  c_n e^{\alpha nu_0(z)} \quad \text{for all} \quad  z\in K.$$ 
Then we have the following estimate of probability for every $\delta>0$
$$\Prob\Big\{Q\in \Pc^n_0, \quad \dist(\mu_Q,\mu_0) \geq \delta\Big\} \leq A n^A c_n e^{-A^{-1}\delta n},$$
where $A=A(K,\mu_0,\alpha)>0$ is a constant independent of $n, c_n$ and $\delta$. 
\end{theorem}

The distances $\dist(\cdot,\cdot)$ and $\dist_\gamma(\cdot,\cdot)$ between probability measures in $\C$ will be introduced in Section \ref{section_abstract}. 
For instance, when $\dist(\mu_Q,\mu_0) \leq \delta$, we have $\dist_\gamma (\mu_Q,\mu_0)\leq c_\gamma \delta^{\gamma/2}$ for $0<\gamma\leq 2$, where $c_\gamma>0$ is a constant, see Proposition \ref{prop_dist_gamma} below. The distance $\dist_1(\cdot,\cdot)$ is equivalent to the well-known Kantorovich-Wasserstein distance. We can deduce from such inequalities estimates of the number of zeros of $Q$ in an open set, see Propositions \ref{prop_sup_mass} and \ref{prop_mass_out_support} below for details. 

Note that in the theorem we can choose the number $\delta$ 
depending on $n$ and the statement
is only interesting when $\log c_n = o(n)$ as $n\to\infty$. We will see later that the case where $c_n$ is a polynomial in $n$ is already interesting for applications in random matrices theory. The above theorem is easier to use when $\alpha$ is an even integer because in this case $|Q(z)|^\alpha$ is a polynomial in $z$ and $\overline z$. 

Note also that in the above theorem the zeros of the polynomial $Q$ do not necessarily lie in the curve $K$. 
Under some extra conditions, we have a more precise estimate of the number of zeros of $Q$ as shown in the following result.

\begin{theorem} \label{th_abstract_bis}
Let $\mu_0, K, u_0, c_n,\alpha$ and the probability measure on $\Pc_0^n$ be as in Theorem \ref{th_abstract}. 
Let $l$ be a smooth curve in $\C$, containing $K$, which is a closed subset of $\C$. 
Let $L$ be a closed subset of $l$ such that in a neighborhood of $L$,
$\mu_0$ has a bounded density with respect to the Lebesgue measure on $l$. 
Assume that the zeros of almost every polynomial $Q\in\Pc_0^n$ lie in $l$. 
For $\delta>0$, denote by $\Fc_n(\delta)$ the set of polynomials $Q\in \Pc_0^n$ such that 
$|\mu_Q(I)-\mu_0(I)| \geq  \sqrt{\delta}$  
for at least one arc  $I\subset L$.
Then we have the following estimate of probability 
$$\Prob(Q\in \Fc_n(\delta)) \leq A n^A c_n e^{-A^{-1}\delta n} \text{ for every } \delta>0,$$
where $A=A(K,l,L,\mu_0,\alpha)>0$ is a constant independent of $n$, $c_n$ and $\delta$.  
\end{theorem}

Note that $n\mu_Q(I)$ is the number of zeros of $Q$ in the interval $I$. As above, the number $\delta$ can depend in $n$ and our estimate is uniform on $I$. 
In the case of Hermitian random matrices, we will choose $l=L=\R$, the real line  in $\C$, since the eigenvalues of a Hermitian matrix are always real. 
We also choose $K=[-2,2]$, the support of the semi-circle distribution $\mu_\scl$ and $\alpha=2$. The distribution $\mu_\scl$ is the probability measure 
supported by $[-2,2]$ with  the density function
$$\rho_\scl(x):={1\over 2\pi} \sqrt{(4-x^2)_+},$$
where $(4-x^2)_+:=\max(4-x^2,0)$.  Its logarithmic potential is given in Lemma \ref{lemma_pot_sc} below.
In order to get precise statements, we need to introduce some notations.

Let $M_n=(\xi_{ij})_{1\leq i,j\leq n}$ be an $n\times n$ Hermitian matrix such that the 
entries $\xi_{ij}$, $1\leq i\leq j\leq n$, are independent random variables. So the diagonal entries $\xi_{ii}$ are real-valued, the off-diagonal ones $\xi_{ij}$ are complex-valued and we have $\xi_{ij}=\overline\xi_{ji}$. 
We will assume that 
\begin{enumerate}
\item $\xi_{ii}$ are of mean 0 and variance bounded by a constant $\alpha\geq 0$;
\item $\xi_{ij}$, with $i\not=j$, are of mean 0, variance 1, and their fourth moments $\E|\xi_{ij}|^4$ are bounded by a constant $\beta>0$.
\end{enumerate}
Define $W_n:=n^{-1/2}M_n$.  
Such a random matrix is called {\it Wigner Hermitian matrix}. 
The {\it empirical spectral distribution (ESD)} of $W_n$ is the probability measure
$$\mu_{W_n}:={1\over n} \sum_{p=1}^n \delta_{\lambda_p},$$
where $\lambda_p$'s denote the eigenvalues of $W_n$. We consider it as a probability measure on $\C$ but it is supported by $\R$ since $W_n$ is Hermitian. It is well-known that when $\alpha$ and $\beta$ are fixed or not too big, $\mu_{W_n}$ converges to $\mu_\scl$, almost surely, as $n$ goes to infinity, see \cite{AGZ, BS, Pastur, Wigner} for more general statements. 
Here, we obtain the following result. 

\begin{theorem} \label{th_main}
With the above notations and hypotheses, there is a universal constant $A>0$  satisfying the following property. For every $\delta>0$,
there is a set  $\Ec_n(\delta)$ of $n\times n$ Hermitian matrices satisfying the following estimate of probability
$$\Prob(M_n\in \Ec_n(\delta))\leq An^A e^{\alpha+\beta/2} e^{-A^{-1}\delta n}$$
and such that if  $M_n\not\in \Ec_n(\delta)$ and $I\subset \R$ is an interval then
$$\dist(\mu_{W_n},\mu_\scl)\leq \delta \qquad \text{and} \qquad |\mu_{W_n}(I)-\mu_\scl(I)| \leq  \sqrt{\delta}.$$
\end{theorem}

Note that our result only requires  a control of the first four moments of the entries of $M_n$. 
When $I$ is outside $[-2,2]$ or is close to $\pm 2$ we can improve the last estimate, see Propositions \ref{prop_mass_out_support} and  \ref{prop_endpoint} below. Note also that the estimates are only significant when $\delta\gg (\log n)n^{-1}$ and the length of $I$ is larger than $\sqrt{\delta}$. 

It is known under stronger hypotheses on the entries of $M_n$ that the smallest size of $I$ for which one has a significant estimate on $\mu(I)$ is $n^{-1}$ up to a logarithmic factor. Many deep works have been done in this direction. We refer the reader to the papers by Erd\"os, Guionnet,  Knowles,  P\'ech\'e, Ram\'irez,  Schlein, Tao, Vu, Yau,  Yin, Zeitouni \cite{Erdos1}-\cite{GZ}, \cite{TaoVu1}-\cite{TaoVu5} and the references therein for details. 

As mentioned above, the almost sure converge to the semi-circular distribution requires a weaker condition on the moments of the entries of $M_n$, see Chapter 2 in Bai-Silverstein \cite{BS}. Theorem \ref{th_main} above gives an alternative for the last step of the proof in this book. 

We will give now an application of the last theorem to the case of matrices associated to random graphs. 
Consider a non-oriented graph with $n$ vertices. Its adjacency matrix is a square $n\times n$ matrix whose diagonal entries are 0 and non-diagonal entry of indexes $(i,j)$ is 1 if there is an edge joining the vertices $i, j$ and 0 otherwise. 
Let $0\leq p\leq 1$ be a real number. Recall that {\it the Erd\"os-R\'enyi random graph} $G(n,p)$ is a graph with $n$ vertices obtained by drawing randomly and independently edges between each pair of vertices with probability $p$.  We will assume that $p\leq 1/2$ because the case $p\geq 1/2$ can be reduced to the first case by considering the complementary graph. 

Denote by $A_n$ the adjacency matrix of $G(n,p)$. It is known that when $np\to\infty$, the empirical spectral distribution of $W_n:={1\over \sqrt{n}\sigma} A_n$ converges in distribution to the semi-circle distribution, where $\sigma:=\sqrt{p(1-p)}$, see \cite{TVW}.  The condition $np\to\infty$ is necessary.  Some explicit estimates on the rate of convergence and the probability to have good rate of convergence were obtained by Tran-Vu-Wang \cite{TVW} and 
Erd\"os-Knowles-Yau-Yin in \cite{Erdos0}. The following corollary is complementary to their results. 

\begin{corollary} \label{cor_graph}
Let $G(n,p), A_n, \sigma$ and $W_n$  be as above with $np\to\infty$ as $n\to\infty$. Let $\mu_{W_n}$ be the empirical spectral distribution of $W_n$. Then there is a universal constant $A>0$ satisfying the following property. For every $\delta>0$,  there is a set $\Gc_n$ of graphs such that the probability that a graph in the model $G(n,p)$ belongs to $\Gc_n$ is at most equal to $An^A e^{1/p} e^{-A^{-1}\delta n}$, and if $W_n$ is associated to a graph outside $\Gc_n$ then 
$$|\mu_{W_n}(I)-\mu_\scl(I)| \leq  \sqrt{\delta}$$
for every interval $I\subset \R$.
\end{corollary}

Note that as above, $\delta$ may depend on $n$ and when $I$ is outside the interval $[-2,2]$ or is close to the endpoints $\pm 2$, we can improve the last estimate in the corollary.

\medskip

The paper is organized as follows. In Section \ref{section_abstract} we will present some basic properties of positive measures and their logarithmic potentials.  We also give there the proofs of Theorems \ref{th_abstract} and \ref{th_abstract_bis}. Theorem \ref{th_main} and Corollary \ref{cor_graph} will be obtained as consequences. The proofs of these results will be given in Section \ref{section_Hermit}.

\medskip
\noindent
{\bf Acknowledgement.} 
This work was supported by the Start-Up 
Grant R-146-000-204-133 from National University of  Singapore.
It was partially written during my visit  at Paris 11 University. I would like to thank Viet-Anh Nguyen and Nessim Sibony for their hospitality and help.

\section{Large deviation theorem for zeros of polynomials} \label{section_abstract}

In this section, we will give some basic properties of probability measures on the complex plane $\C$, the distances between them and 
their logarithmic potentials. We refer the reader to \cite{Hormander, Tsuji} for the complex potential theory.
We then give the proofs of Theorems \ref{th_abstract} and \ref{th_abstract_bis}.

Recall that the Riemann sphere $\P^1=\C\cup\{\infty\}$ is the natural compactification of $\C$ by adding a point $\infty$ at infinity. 
We will use the standard Hermitian metric on $\P^1$ given by the Fubini-Study form $\omega_\FS:=\ddc  \log (1+|z|^2)^{1/2}$, where $z$ denotes the standard complex coordinate in $\C$ and $\ddc:={i\over \pi} \ddbar$. This differential form extends to a smooth form on $\P^1$.  

Let $\Mc_c(\C)$ denote the set of all probability measures with compact support in $\C$. For each $\mu\in\Mc_c(\C)$, there is a unique subharmonic function $u:\C\to \R\cup\{-\infty\}$ such that 
$$\ddc u=\mu \quad \text{and} \quad \lim_{z\to \infty}u(z)-\log|z|=0.$$
The first identity is understood in the sense of currents or distributions. If $u$ is of class $\Cc^2$ on $\C$, then $\ddc u = \Delta u {i\over 4\pi} dz\wedge d\overline z={1\over 2\pi}\Delta u dx\wedge dy$ which is identified to a measure with density ${1\over 2\pi} \Delta u$ on $\C$. Here, we write $z=x+iy$ with $x,y\in\R$. The function $u$ is called {\it the logarithmic potential} of $\mu$. It is not difficult to see that  $u$ is given by the formula
$$u(z):=\int_\C \log|z-a| d\mu(a).$$

Let $\mu,\mu'$ be two probability measures in $\Mc_c(\C)$ and $u,u'$ their logarithmic potentials. Define
$$\dist(\mu,\mu'):=\|u-u'\|_{L^1(\P^1)}:=\int_{\P^1} |u-u'| \omega_\FS.$$
It is not difficult to check that $\dist(\cdot,\cdot)$ is a distance on  $\Mc_c(\C)$. 
For every $\gamma>0$, define
$$\dist_\gamma(\mu,\mu'):=\sup \Big\{ |\langle \mu-\mu',\phi\rangle|, \ \phi \text{ is a } \Cc^\gamma \text{ function on } \P^1 \text{ with } \|\phi\|_{\Cc^\gamma}\leq 1\Big\},$$
where the pairing $\langle \mu-\mu',\phi\rangle$ denotes the integral of $\phi$ with respect to the measure $\mu-\mu'$. Recall from \cite{DinhSibony4} that for $0<\gamma\leq \gamma'$, there is a constant $c>0$ depending only on $\gamma$ and $\gamma'$ such that 
$$\dist_{\gamma'}\leq \dist_\gamma\leq c[\dist_{\gamma'}]^{\gamma/\gamma'}.$$
The function $\dist_\gamma$ is also a distance  on $\Mc_c(\C)$. Note that $\dist_1$ is equivalent to the Kantorovich-Wasserstein distance. 
We have the following proposition.

\begin{proposition} \label{prop_dist_gamma}
For every $0<\gamma\leq 2$, there is a constant $c_\gamma>0$ such that 
$$\dist_\gamma\leq c_\gamma \dist^{\gamma/2}.$$
\end{proposition}
\proof
Using the above comparison between $\dist_\gamma$ and $\dist_{\gamma'}$, we reduce the problem to the case where $\gamma=2$. Let $\phi$ be any $\Cc^2$ function on $\P^1$ such that $\|\phi\|_{\Cc^2}\leq 1$. With the notations as above,  using Stokes formula on $\P^1$, we have 
$$|\langle \mu- \mu',\phi\rangle| = |\langle \ddc (u-u'),\phi\rangle |= |\langle u-u',\ddc \phi\rangle|.$$

Since  $\|\phi\|_{\Cc^2}\leq 1$, we can write $\ddc\phi=\omega^+-\omega^-$ where $\omega^\pm$ are positive $(1,1)$-forms on $\P^1$, bounded by a positive constant $c$ times $\omega_\FS$. We then deduce from the last identities that 
$$|\langle \mu- \mu',\phi\rangle| \leq 2c \|u-u'\|_{L^1(\P^1)}= 2c \dist(\mu,\mu').$$
The proposition follows.
\endproof

\begin{proposition} \label{prop_dist_sup}
Let $\mu_0$ be a probability measure whose support is contained in a compact smooth non-closed curve $K$ in $\C$. Let  $u_0$ be the logarithmic potential of $\mu_0$. 
Assume that $u_0$ is continuous. There is a constant $c>0$ depending only on $K$ and $\mu_0$ satisfying the following property.
If $\mu$ is a probability measure with compact support in $\C$ and $u$ is its logarithmic potential,  then
$$\dist(\mu,\mu_0) \leq c\sup_K (u-u_0)=c\sup_{\P^1}(u-u_0).$$ 
In particular, for every $0<\gamma\leq 2$, there is a constant $c_\gamma>0$ depending only on $K,\mu_0$ and $\gamma$ such that 
$$\dist_\gamma(\mu,\mu_0)\leq c_\gamma\sup_K (u-u_0)^{\gamma/2}=c_\gamma\sup_{\P^1} (u-u_0)^{\gamma/2}.$$ 
\end{proposition}
\proof
By Proposition \ref{prop_dist_gamma}, we only need to prove the first assertion. Define 
$$\delta:=\sup_K (u-u_0) \qquad \text{and} \qquad v:=u-u_0-\delta.$$ 
Recall that since $\mu_0$ is supported by $K$, the function $u_0$ is harmonic in $\C\setminus K$. 
So $v$ is negative on $K$ and subharmonic on $\P^1\setminus K$ with $v(\infty)=-\delta$. By maximum principle, 
we deduce that $v$ is negative on $\P^1$ and therefore
$$\sup_K (u-u_0)=\sup_{\P^1}(u-u_0)=\delta.$$ 
We also have 
$$\dist(\mu,\mu_0)=\int_{\P^1} |u-u_0| \omega_\FS \leq \Big|\int_{\P^1} \big(\delta + |\delta - (u-u_0)|\big) \omega_\FS\Big|=\delta  -\int_{\P^1} v\omega_\FS.$$
Hence, we only need to bound the last integral from below.
 
Consider the conformal bijective map $\pi:\D\to\P^1\setminus K$ such that $\pi(0)=\infty$, where $\D$ is the unit disc in $\C$. 
This map is Lipschitz up to the boundary, see \cite[Th.3.9]{Pommerenke}. Therefore, the form $\pi^*(\omega_\FS)$ is bounded on $\D$. Define $\hat v:= v\circ\pi$. We have 
$$\int_{\P^1} v\omega_\FS =\int_\D \hat v \pi^*(\omega_\FS) \geq c\int_\D \hat v d\Leb,$$
where $c>0$ is a constant and $\Leb$ denotes the Lebesgue measure. Now, since $\hat v$ is a subharmonic function with $\hat v(0)=v(\infty)=-\delta$, the submean inequality implies that the last integral is bounded from below by  $-\pi \delta$. The proposition  follows.
\endproof

If $I$ is a subset of $\C$ and $\epsilon>0$, denote by $I_\epsilon$ the open set of points $z\in \C$ such that $\dist(z,I)<\epsilon$. We have the following results.

\begin{lemma} \label{lemma_sup_L1}
Under the hypotheses of Proposition  \ref{prop_dist_sup}, define $\epsilon:=\sqrt{\sup_K (u-u_0)}$.  Then, if $\epsilon\leq 1$ we have 
$$\int_{K_\epsilon}|u(z)-u_0(z)| d\Leb(z) \leq c\epsilon^3,$$
where $c>0$ is a constant depending only on $\mu_0$ and $K$. 
\end{lemma}
\proof
Observe that $u-u_0$ is subharmonic on $\P^1\setminus K$ and vanishes at $\infty$.  We have seen in Proposition \ref{prop_dist_sup} that  $u\leq u_0+\epsilon^2$ on $\P^1$. So $u$ belongs to a compact family of subharmonic functions on $\C$.
Here, the compactness is with respect to the $L^1_\loc$ topology. Therefore, in order to prove the lemma, we only need to consider $\epsilon$ small enough.

Choose a smooth compact closed curve $K'$ containing $K$. Let $\Omega$ be the unbounded connected component of $\C\setminus K'$ and define $K_\epsilon^+:=K_\epsilon\cap \Omega$. Let $\pi:\D\to \Omega\cup\{\infty\}$ be a bijective conformal map such that $\pi(0)=\infty$. 
By Kellogg-Warschawski theorem, this map is smooth up to the boundary, see \cite[Th.3.6]{Pommerenke}. Then, there is a constant $c'>0$ such that $\pi^{-1}(K_\epsilon^+)$ is contained in the annulus $A_\epsilon:=\{1-c' \epsilon\leq |z|\leq 1\}$. It follows that there is a constant $c''>0$ such that if we define $v:=u-u_0-\epsilon^2$ and $\tilde v:= v\circ\pi$ then using that $v\leq 0$ we have
\begin{eqnarray*}
\int_{K_\epsilon^+}|u(z)-u_0(z)| d\Leb(z) & = & \int_{K_\epsilon^+}|\epsilon^2+v(z)| d\Leb(z) \ \leq\  \int_{K_\epsilon^+} (\epsilon^2-v(z) ) d\Leb(z) \\
&  \leq &  c''\epsilon^3 - c''\int_{A_\epsilon} \tilde v(z) d\Leb(z).
\end{eqnarray*}

On the other hand, since the function $\tilde v$ is subharmonic, its average in each circle of center 0 is larger than or equal to $\tilde v(0)=-\epsilon^2$. So the last integral is bounded from below by a negative constant times $\epsilon^3$. It follows that 
$$\int_{K_\epsilon^+}|u(z)-u_0(z)| d\Leb(z)\leq {1\over 2} c\epsilon^3,$$
for some constant $c>0$. 

To complete the proof, we need a similar estimate for $K_\epsilon^-:=K_\epsilon\setminus K_\epsilon^+$. Choose another smooth closed curve $K''$ equal to $K'$ in a neighborhood of $K$ such that $K_\epsilon^-$ is exterior to $K''$. Arguing as above, we obtain the desired estimate for $K_\epsilon^-$.
\endproof

Recall that we don't assume that $\mu$ is supported by $K$.

\begin{proposition} \label{prop_sup_mass}
Let $K,\mu_0,u_0, \mu, u$ be as in Proposition \ref{prop_dist_sup}. Define $\epsilon:=\sqrt{\sup_K (u-u_0)}$. Then there is a constant $c>0$ depending only on $K$ and $\mu_0$ such that for any arc $I\subset K$ we have
$$\mu_0(I)-c\epsilon\leq \mu(I_\epsilon)\leq \mu_0(I_{2\epsilon}\cap K) + c \epsilon.$$
In particular, we have $\mu(\C\setminus K_\epsilon)\leq c\epsilon$. 
\end{proposition}
\proof
Since $\mu$ and $\mu_0$ are probability measures and $\mu_0(K)=1$, the second assertion is a direct consequence of the first inequality applied to $K$ instead of $I$. 
We prove now the first assertion. 
Observe that we only need to consider $\epsilon$ small enough. 
We have seen in Proposition \ref{prop_dist_sup} that $u\leq u_0+\epsilon^2$ on $\P^1$. 
Therefore, if $L$ is a smooth compact curve containing $K$, we have 
$\sup_L(u-u_0)=\epsilon^2$. Replacing $K$ by a larger curve allows us to assume that the endpoints of $K$ are outside $I_{2\epsilon}$.

Let $0\leq \chi\leq 1$ be a smooth cut-off function equal to 1 on $I_\epsilon$ and to 0 outside $I_{2\epsilon}$. We can choose $\chi$ so that $\|\chi\|_{\Cc^2}$ is bounded by a constant times $\epsilon^{-2}$. We have
$$\mu(I_\epsilon)\leq \langle \mu, \chi\rangle =\langle \mu_0,\chi\rangle +\langle \mu-\mu_0,\chi\rangle\leq \mu_0(I_{2\epsilon}\cap K)+ \langle \mu-\mu_0,\chi\rangle.$$
On the other hand, by Stokes formula and using that $\chi$ has support in $K_{2\epsilon}$, we have
$$\langle \mu-\mu_0,\chi\rangle=\langle \ddc (u-u_0),\chi\rangle = \langle u-u_0,\ddc \chi\rangle \leq \|\chi\|_{\Cc^2} \int_{K_{2\epsilon}} |u-u_0| d\Leb.$$
The estimate on the $\Cc^2$ norm of $\chi$ and Lemma \ref{lemma_sup_L1} imply that $\langle \mu-\mu_0,\chi\rangle=O(\epsilon)$. It follows that 
$$\mu(I_\epsilon) \leq  \mu_0(I_{2\epsilon}\cap K) + c\epsilon$$
for some constant $c>0$. 

To complete the proof of the first assertion, we need to check that 
$\mu(I_\epsilon)\geq \mu_0(I)-c\epsilon$ for some constant $c>0$.
Let $0\leq \chi'\leq 1$ be a cut-off function equal to 1 on $I$ and 0 outside $I_\epsilon$. We can choose such a function with $\Cc^2$ norm bounded by a constant times $\epsilon^{-2}$. As above, we obtain
$$\mu(I_\epsilon)\geq \langle \mu,\chi'\rangle \geq \mu_0(I)- \|\chi'\|_{\Cc^2} \int_{K_{\epsilon}} |u-u_0| d\Leb=\mu_0(I)-O(\epsilon).$$
The proposition follows.
\endproof

\begin{corollary} \label{cor_sup_mass}
Let $K,\mu_0,u_0, \mu, u, \epsilon$ be as in Proposition \ref{prop_sup_mass}.  Let $l$ be a smooth curve in $\C$ containing $K$ which is a closed subset of $\C$. Let $L$ be a closed subset of $l$ such that in an open neighborhood of $L$, $\mu_0$ has bounded density with respect to the Lebesgue measure on $l$.
Let $U$ be an open neighborhood of $K$ and assume that the support of $\mu$ is contained in $l\cup (\C\setminus U)$. Then for any arc $I\subset L$ we have 
$$|\mu(I)-\mu_0(I)| \leq c\epsilon,$$
where $c=c(K,U,\mu_0,l,L)>0$ is a constant independent of $\mu, I$ and $\epsilon.$
\end{corollary}
\proof
We only need to consider small $\epsilon$. Choose an open neighborhood $V$ of $K$ such that $V\Subset U$ and the intersection $K':=l\cap \overline V$ is connected.
By the last assertion of Proposition \ref{prop_sup_mass} or Proposition \ref{prop_mass_out_support} below, we can replace $I$ by $I\cap \overline V$ in order to assume that $I\subset K'\cap L$. 

We have seen in  Proposition \ref{prop_sup_mass}  that $\epsilon=\sqrt{\max_{K'} (u-u_0)}$. 
So we can apply Proposition \ref{prop_sup_mass} to $K'$ instead of $K$. Using that the density of $\mu_0$ is bounded in a neighborhood of $L$ (which contains $I_{2\epsilon}$ because $\epsilon$ is small), we have 
$$\mu(I)\leq \mu(I_\epsilon) \leq \mu_0(I_{2\epsilon}\cap K')+O(\epsilon)=\mu_0(I)+O(\epsilon).$$ 
It remains to check that $\mu(I)\geq \mu_0(I)-O(\epsilon)$. Observe that it is enough to consider the case where the length of $I$ is at least $3\epsilon$ because otherwise $\mu(I)\geq 0 \geq \mu_0(I)-O(\epsilon)$. 

If $a$ and $b$ are the endpoints of $I$, denote by $I'$ the set of points $z\in I$ such that $\dist(z,a)\geq \epsilon$ and $\dist(z,b)\geq \epsilon$. 
Applying Proposition \ref{prop_sup_mass} to $I'$ and using that $\mu$ has support in $l\cup (\C\setminus U)$, we get
$$\mu(I)=\mu(I'_\epsilon) \geq \mu_0(I')-O(\epsilon) =\mu_0(I)-O(\epsilon).$$
The corollary follows.
\endproof

\begin{remark}\rm  \label{rk_density}
Note that when $\mu_0$ is not of bounded density,  we can work with suitable neighborhood of $I$ instead of using $I_\epsilon$, in order to get better estimates. The choice of this neighborhood  depends on the measure $\mu_0$. Similarly, we can also improve the estimates on the arcs where the density of $\mu_0$ is small, as shown in  Proposition \ref{prop_endpoint} below.  
\end{remark}

The following result is complementary to the last estimate in Proposition \ref{prop_sup_mass}.

\begin{proposition} \label{prop_mass_out_support}
Let $K,\mu_0,u_0, \mu, u$ be as in Proposition \ref{prop_sup_mass}. Let $U$ be a neighborhood of $K$. Then there is a constant $c>0$ depending only on $K,\mu_0$ and $U$ such that  
$$\mu(\C\setminus U) \leq c\sup_K (u-u_0).$$
\end{proposition}
\proof
Fix a smooth cut-off function $0\leq \chi\leq 1$ with compact support in $U$ which is equal to 1 on $K$. We have 
$$\mu(U) \geq \langle \mu,\chi\rangle =\langle \mu_0,\chi\rangle +\langle \mu-\mu_0,\chi\rangle=1 +\langle u- u_0,\ddc \chi\rangle.$$
By the first inequality in Proposition \ref{prop_dist_sup}, the last integral is larger than a negative constant times  $\sup_K (u-u_0)$. We then obtain the result using that $\mu(\C\setminus U)=1-\mu(U)$ because $\mu$ is a probability measure.
\endproof

In order to prove Theorems \ref{th_abstract} and \ref{th_abstract_bis}, 
we will need the following version of Markov brother's inequality, due to Szeg\"o, see also \cite[p.567]{Schaeffer}. Note that instead of this inequality we can also use a version of Bernstein-Markov inequality to get a slightly weaker estimate but valid in a more general setting, see \cite[Cor. 3.13]{DinhNguyen} for details.

\begin{lemma} \label{lemma_markov_0}
Let $K$ be the graph in $\C\simeq \R^2$, of a $\Cc^2$ real-valued function defined on $[0,1]$ whose $\Cc^2$ norm is bounded by a fixed constant $M$. Then there is a constant $c>0$ depending only on $M$ such that if $Q$ is any polynomial of degree $n\geq 1$ we have 
$$\max_K |Q'| \leq cn^2 \max_K|Q|.$$ 
\end{lemma}
\proof
We give here the proof for the reader's convenience, see  \cite[p.567]{Schaeffer}. 
Multiplying $Q$ by a constant allows us to assume that $\max_K |Q|=1$. 
Let $\phi:\P^1\setminus K\to \P^1\setminus \overline\D$ be a conformal bijective map  such that $\phi(\infty)=\infty$. 
So $Q(z)\phi(z)^{-n}$ is holomorphic on $\P^1\setminus K$ and has modulus bounded by 1 when $z\to K$. By maximum modulus principle, we have 
$|Q(z)||\phi(z)|^{-n}\leq 1$ on $\P^1\setminus K$ and hence $|Q(z)|\leq |\phi(z)|^n$ on $\P^1\setminus K$. Observe that
there is a constant $\gamma>0$ such that $|\phi(z)|-1\leq \gamma \dist (z,K)^{1/2}$ for $z$ is a fixed large compact subset of $\C$, see \cite[Th.3.9]{Pommerenke}.  We then deduce from Cauchy's type formula that for every $z_0\in K$
\begin{eqnarray*}
|Q'(z_0)| & = & {1\over 2\pi} \Big|\int_{|z-z_0|  =\gamma^{-2}n^{-2}} {Q(z)dz\over (z-z_0)^2}\Big| \ \leq \  \gamma^{2}n^2
 \max_{|z-z_0| =   \gamma^{-2}n^{-2}} |Q(z)| \\
 &  \leq &  \gamma^2n^2 \max_{|z-z_0|= \gamma^{-2}n^{-2}} |\phi(z)|^n. 
\end{eqnarray*}
For $z$ as in the last line, we have 
$$|\phi(z)|^n\leq  (1+ \gamma\dist(z,K)^{1/2})^n \leq  (1+ \gamma |z-z_0|^{1/2})^n = (1+n^{-1})^n\leq e.$$
The result follows.
\endproof

\begin{lemma} \label{lemma_markov}
Let $K$ be a compact smooth curve in $\C$ and $\phi$  a real-valued H\"older continuous function on $K$. There is a constant $A>0$ satisfying the following property. If $Q$ is any polynomial of degree $n\geq 1$ and if  $\Sigma\subset K$ is an arbitrary subset such that $\dist(z,\Sigma)\leq A^{-1}n^{-A}$ for every $z\in K$, then we have
$$\max_{z\in K} |Q(z)| e^{n\phi(z)} \leq 2 \max_{z\in \Sigma} |Q(z)| e^{n\phi(z)}.$$
\end{lemma}
\proof
Fix a constant $A>2$ large enough so that if $z,z'\in K$ and $|z-z'|\leq A^{-1/3} n^{-A+2}$, we have $|\phi(z)-\phi(z')|\leq 1/(4n)$. 
Choose $a\in K$ such that 
$$ |Q(a)| e^{n\phi(a)}=\max_{z\in K} |Q(z)| e^{n\phi(z)}.$$
Choose a compact arc $L\subset K$ of diameter $A^{-1/3}n^{-A+2}$ containing $a$. 

The above inequality on $\phi$ implies that for $z\in L$
$$|\phi(z)-\phi(a)|\leq{1\over 4n} \qquad \text{and hence} \qquad  |Q(z)| e^{n\phi(z)} \geq {1\over \sqrt{2}}  |Q(z)| e^{n\phi(a)}.$$
So we only need to check that 
$$|Q(a)|\leq \sqrt{2} \max_{z\in\Sigma\cap L}  |Q(z)|.$$

Denote by $\tilde L$ and $\tilde \Sigma$  the images of $L$ and $\Sigma$ by a suitable affine map $\pi(z)= \alpha z+\beta$ so that $\tilde L$ is the graph over $[0,1]$ of a function with bounded $\Cc^1$ norm and vanishing at $0$ and $1$.  We have $|\alpha|=A^{1/3}n^{A-2}$. 
Define also
$\tilde Q(z):= Q(\pi^{-1}(z))$. It is enough to check that 
$$\max_{z\in \tilde L} |\tilde Q(z)| \leq \sqrt{2}  \max_{z\in\tilde\Sigma\cap \tilde L}  |\tilde Q(z)|.$$
Let $b\in \tilde L$ such that $|\tilde Q(z)|$, $z\in \tilde L$, archives its maximum at $b$. By hypothesis on $\Sigma$, there is a point $b'\in \tilde\Sigma\cap \tilde L$ such that the arc in $\tilde L$ joining $b$ and $b'$ has length smaller than $A^{-1/2}n^{-2}$. 
Denote by $\gamma(b,b')$ this length.
By Lemma \ref{lemma_markov_0}, we have $$\max_{z\in \tilde L}|\tilde Q'(z)|\leq cn^2\max_{z\in \tilde L}|\tilde Q(z)| = cn^2 |\tilde Q(b)|.$$
Since $A$ is large enough, we easily deduce that 
$$|\tilde Q(b)|\leq |\tilde Q(b')|+cn^2 |\tilde Q(b)|\gamma(b,b') \leq |\tilde Q(b')|+ (1-1/\sqrt{2})|\tilde Q(b)|.$$
Hence, $|\tilde Q(b)|\leq \sqrt{2} |\tilde Q(b')|$. The lemma follows.
\endproof

\noindent
{\bf Proof of Theorem \ref{th_abstract}.} Fix a constant  $A'>0$ large enough.
Let $\Sigma$ be a set of about $A'n^{A'}$ points equidistributed on $K$. 
Consider the set $\Fc$ of polynomials $Q\in\Pc^n_0$ such that 
$$|Q(z)|^\alpha\geq 2^{-\alpha} e^{A'^{-1}\alpha n\delta} e^{\alpha n u_0(z)}$$
for at least one point $z\in \Sigma$. By hypothesis on the expectation of $|Q|^\alpha$, the size of $\Fc$ satisfies
$$\Prob(Q\in \Fc) \leq 2^\alpha A'n^{A'}c_ne^{-A'^{-1}\alpha \delta n}.$$

Consider now an arbitrary $Q\not\in \Fc$. We have $|Q(z)| \leq {1\over 2}e^{A'^{-1}n\delta} e^{nu_0(z)}$ for $z\in \Sigma$. 
By Lemma \ref{lemma_markov}, applied to the function $-u_0(z)$ instead of $\phi(z)$, we obtain that 
$$|Q(z)| \leq e^{A'^{-1}n\delta} e^{nu_0(z)} \quad \text{for}\quad  z\in K.$$ 
So the logarithmic potential $u$ of the measure $\mu_Q$ satisfies
$$u={1\over n} \log |Q(z)| \leq u_0(z) + A'^{-1}\delta  \quad \text{for}\quad  z\in K.$$ 
By Proposition \ref{prop_dist_sup}, we have $\dist(\mu_Q,\mu_0)<\delta$ since $A'$ is large enough. Finally,  by taking a constant $A$ large enough, we obtain the result.
\hfill $\square$

\medskip

Note that we can also apply Propositions \ref{prop_sup_mass} and \ref{prop_mass_out_support} to $Q\not\in\Fc$ as above.

\medskip

\noindent
{\bf Proof of Theorem \ref{th_abstract_bis}.} 
Consider $\Fc$ and $Q\not\in\Fc$ as in the proof of Theorem \ref{th_abstract}. 
It is then enough to apply Corollary \ref{cor_sup_mass} to the function $u:={1\over n} \log|Q|$ and to $\epsilon:=\sqrt{A'^{-1}\delta}$. 
\hfill $\square$

\medskip

As mentioned in Remark \ref{rk_density}, we can improve the estimate near the points where the density of $\mu_0$ is small. For simplicity, we consider the case of the semi-circle distribution. We have the following proposition.

\begin{proposition} \label{prop_endpoint}
Under the hypotheses of Corollary \ref{cor_sup_mass}, assume for simplicity that $\mu_0=\mu_\scl$, $K=[-2,2]$, $l=L=\R$ and $U=\C$. Then there is a universal constant $c>0$ such that 
$$\mu\big(\R\setminus [-2+\epsilon^{4/5}, 2-\epsilon^{4/5}]\big) \leq c \epsilon^{6/5}.$$ 
\end{proposition}
\proof
We can assume that $\epsilon$ is small. 
It is enough to show that $\mu((-\infty,-2+\epsilon^{4/5})) = O(\epsilon^{6/5})$ since we can obtain in the same way that 
 $\mu((2-\epsilon^{4/5}, +\infty)) = O(\epsilon^{6/5})$.
 Using Proposition \ref{prop_mass_out_support}, we only need to show that $\mu([-3,-2+\epsilon^{4/5}]) = O(\epsilon^{6/5})$.
We will use that the density of $\mu_\scl$ is small near $\pm 2$. 

Define $\eta:=\epsilon^{4/5}$ and $I:=[-3,-2+\eta]$. As in Proposition \ref{prop_sup_mass}, using a cut-off function in $I_\eta$, we obtain that 
$$\mu(I)\leq \mu_0(I_\eta) + O(\eta^{-2}) \int_{I_\eta} |u-u_0| d\Leb.$$
The explicit form of $\mu_\scl$ implies that the first term in the last sum is $O(\eta^{3/2})=\epsilon^{6/5}$ and as in Lemma \ref{lemma_sup_L1}, we obtain that the second term of the last sum is $O(\eta^{-2}) \eta\epsilon^2=\epsilon^{6/5}$. 
This implies the result. 
\endproof

\section{Local semi-circle law for Hermitian random matrices} \label{section_Hermit}

In this section, we give the proofs of Theorem \ref{th_main} and Corollary \ref{cor_graph} as a consequence of the abstract results obtained in Theorems \ref{th_abstract} and \ref{th_abstract_bis} for $\alpha=2$, $K=[-2,2]$ and $\mu_0=\mu_\scl$. We will use the notations introduced above. Let $u_\scl$ denote the logarithmic potential of $\mu_\scl$.
So we need to compare $\E(\det (z-W_n)^2)$ with $e^{2nu_\scl(z)}$ for $z\in [-2,2]$, see Proposition \ref{prop_expectation} below.

Consider $\pi(w)=w+1/w$ which defines  bijective conformal maps from the unit disc $\D\subset\C$ and $\P^1\setminus \overline\D$ to $\P^1\setminus [-2,2]$. We also have $\pi(0)=\pi(\infty)=\infty$ and that $\pi$ is a 2:1 map from the unit circle to $[-2,2]$. Define $\hat u:=u\circ\pi$. The following lemma gives us the explicit values of $u_\scl$.

\begin{lemma} \label{lemma_pot_sc}
We have 
$$\hat u_\scl(w)={1\over 4}( w^2+\overline w^2)-\log|w| \quad \text{for} \quad |w|\leq 1$$
and
$$\hat u_\scl(w)={1\over 4}( w^{-2}+\overline w^{-2})+\log|w| \quad \text{for} \quad |w|\geq 1.$$
In particular, we have 
 $$u_\scl(z)= {1\over 4} (z^2-2) \quad \text{for} \quad z\in [-2,2].$$
\end{lemma}
\proof
The last assertion is a direct consequence of the first one by taking $z=w+1/w$ with $|w|=1$. 
The second identity is a consequence of the first one because $\hat u_\scl(w)= \hat u_\scl(1/w)$. 
We prove now the first identity. 
Write $z=w+1/w$ with $|w|\leq 1$. Then we have
$$\hat u_\scl(w)=u_\scl(z)={1\over 2\pi} \int_{[-2,2]} \log |z-a| \sqrt{4-a^2} da.$$
Write $a=b+1/b$ with $b=e^{i\vartheta}$ and $\vartheta\in[0,\pi]$. Then the last identities imply
\begin{eqnarray*}
\hat u_\scl(w) & = &{2\over \pi} \int_0^\pi\log\big| w+1/w-e^{i\vartheta} - e^{-i\vartheta}\big| \sin^2\vartheta d\vartheta \\
& = & {2\over \pi} \int_0^\pi \big( \log|1-we^{i\vartheta}| +\log |1-we^{-i \vartheta}| -\log |w|\big) \sin^2 \vartheta d\vartheta \\
& = & {2\over \pi}\int_{-\pi}^\pi \log |1-we^{i\vartheta}| \sin^2\vartheta d\vartheta -\log|w|.
\end{eqnarray*}

We need to show that the first term in the last line is equal to ${1\over 4}( w^2+\overline w^2)$. This  term is equal to the real part of
$${2\over \pi}\int_{b\in \S^1} \log (1-wb) \Big({b-1/b\over 2i}\Big)^2 (-ib^{-1} db)=-{1\over 2i\pi}\int_{b\in \S^1} \log (1-wb) (b-2b^{-1}+b^{-3})db,$$
where we consider the principal branch of the complex logarithmic function.
Now, the residue formula with a unique pole at 0 gives
$$-{1\over 2i\pi }\int_{b\in \S^1} \log (1-wb) (b-2b^{-1}+b^{-3})db={1\over 2}w^2.$$
Taking the real part gives the result.
\endproof

We will need the following lemma.

\begin{lemma} \label{lemma_Hermite}
Let $H_n$ be the (physicists') Hermite polynomial of degree $n$. Then, we have 
$$\Big|H_n\big(\sqrt{n\over 2}z\big)\Big|\leq 
c(2n)^{{1\over 2}(n+1)} e^{n({1\over 4}z^2-{1\over 2})} \quad  \text{for} \quad  z\in\R,$$
where $c>0$ is a universal constant. 
\end{lemma}
\proof
Recall that the Hermite polynomial $H_n(z)$ satisfies
$$H_n(z)=\sum_{0\leq k\leq n/2} {(-1)^k n! \over k! (n-2k)!}(2z)^{n-2k}={n!\over 2i\pi} \oint{e^{2tz-t^2}\over t^{n+1} }dt,$$
where the integral is taken on any simple loop about the origin.  
So we have 
$$H_n\big(\sqrt{n\over 2}z\big)={n!\over 2i \pi}\oint {e^{\sqrt{2n} zt-t^2}\over t^{n+1}} dt.$$
By Stirling's approximation, we have for some universal constant $c>0$
$$(2n)^{-(n+1)/2}\Big|H_n\big(\sqrt{n\over 2}z\big)\Big|\leq c\Big({n\over 2}\Big)^{n/2} e^{-n}\Big|\oint {e^{\sqrt{2n} zt-t^2}\over t^{n+1}} dt\Big|
=c\Big|\oint {e^{n(zs-{1\over 2}s^2-1)}\over s^{n+1}} ds\Big|,$$
where we used the change of variable $t=\sqrt{{n\over 2}} s$. 

Now, consider the last integral on the unit circle in $\C$ and write $s=e^{i\theta}$ with $\theta\in [0,2\pi]$. Since $z$ is real, the real part of $zs-{1\over 2}s^2-1$ is equal to $z\cos \theta-\cos^2\theta-{1\over 2}\cdot$ So it is bounded from above by ${1\over 4}z^2-{1\over 2}$. The lemma follows.
\endproof

We continue the proof of Theorem \ref{th_main}. 
Define 
$$P_{M_n}(z) : = \det (z- M_n) \quad \text{and} \quad Q_{W_n}(z):=\det (z-W_n)=n^{-n/2} P_{M_n}(n^{1/2} z).$$
Recall that we want to bound $\E(|Q_{W_n}(z)|^2)$. 
We need to introduce some notations. 

For any square matrix $M$, denote by $M{[0]}$ the matrix obtained from $M$ by replacing the diagonal entries with 0. 
Define $\llbracket1,n\rrbracket:=\{1,\ldots,n\}$. 
Let $I=\{i_1,\ldots,i_m\}$ be any subset of $\llbracket1,n\rrbracket$. 
The integer $m$ is the length of $I$ which will be also denoted by $|I|$. 
If $M$ is a square $n\times n$ matrix, denote by $M[I]$ the minor obtained  from $M$ by deleting the lines and columns of indexes $i\in I$. 

Recall that $M_n=(\xi_{ij})_{1\leq i,j\leq n}$. Denote by $\Sym(I)$ the group of all permutations of $I$ and  define
$$\xi_I:=\prod_{1\leq k\leq m} \xi_{i_ki_k} \quad \text{and}\quad \xi_\sigma:=\prod_{1\leq k\leq m} \xi_{i_k\sigma(i_k)}\quad \text{for}\quad 
\sigma\in \Sym(I).$$
Let  $\Sym_{\geq k}(I)$ (resp. $\Sym_k(I)$) be the sets of $\sigma\in \Sym(I)$ containing no cycle of length less than $k$ (resp. containing only cycles of length $k$). 

We have the following expansion of determinant
$$P_{M_n{[0]}}(z)=\sum_{\sigma\in\Sym_{\geq 2}(I)} (-1)^{|I|+\sign(\sigma)}\xi_\sigma z^{n-|I|},$$
where $\sign(\sigma)$ denotes the signature of $\sigma$.
Note that the above sum is taken over $\sigma$ and also over the subset $I$ of $\llbracket 1,n\rrbracket$. 

\begin{lemma} \label{lemma_exp_M}
The expectation of $P_{M_n{[0]}}(z)$ satisfies
$$\E\big(P_{M_n{[0]}}(z)\big)=2^{-n/2} H_n\Big({z\over \sqrt{2}}\Big).$$
In particular, there is a universal constant $c>0$ such that we have 
$$\big|\E\big(P_{M_n{[0]}}(z)\big)\big|^2 \leq c n^{n+1} e^{{1\over 2} z^2-n} \quad \text{for} \quad z\in\R.$$
\end{lemma}
\proof
In the expansion of $P_{M_n[0]}(z)$ if a term contains $\xi_{ij}$ but not $\xi_{ji}$, its expectation vanishes because $\xi_{ij}$ has zero mean. We then deduce that the expectation of $P_{M_n{[0]}}(z)$ is equal to the one of 
$$\sum_{\sigma\in\Sym_{2}(I)} (-1)^{|I|+\sign(\sigma)}\xi_\sigma z^{n-|I|}.$$
If $|I|$ is odd, the set $\Sym_{2}(I)$ is empty. So we only need to consider $I$ with even length.

Denote by $2l$ the length of $I$.  
Then we have $\sign(\sigma)=l$ mod 2 for $\sigma\in\Sym_2(I)$. The number of sets $I$ of length $2l$ is ${n\choose 2l}$. 
For such a set  $I$, the cardinality of $\Sym_2(I)$ is $(2l-1)!!$. 
The expectation of $\xi_\sigma$ is 1 because $\xi_{ij}=\overline\xi_{ji}$ and the variance of $\xi_{ij}$ is 1 for $i\not=j$. It follows that the expectation to compute is equal to
$$\sum_{l\leq n/2} (-1)^l {n\choose 2l} (2l-1)!! z^{n-2l}=\sum_{l\leq n/2} {(-1)^l n!\over 2^l l! (n-2l)!} z^{n-2l}=2^{-n/2} H_n\Big({z\over \sqrt{2}}\Big).$$
This gives the first assertion in the lemma. The second one is then deduced from Lemma \ref{lemma_Hermite} above.
\endproof

Consider the expansion
$$|P_{M_n{[0]}}(z)|^2=\sum_{\sigma\in\Sym_{\geq 2}(I)\atop \sigma'\in\Sym_{\geq 2}(I')} (-1)^{|I|+|I'|+\sign(\sigma)+\sign(\sigma')}\xi_\sigma \xi_{\sigma'} z^{n-|I|} \overline z^{n-|I'|}.$$
We will consider some partial sums. 
Denote by $R^{[0]}_{M_n}(z)$ the sum of all terms in the last expression with $\sigma\in\Sym_2(I)$ and $\sigma'\in\Sym_2(I')$ having no common 2-cycle. 
Denote by $R^{[1]}_{M_n}(z)$ the sum of the terms  with $\sigma\in\Sym_{\geq 2}(I)$ and $\sigma'\in\Sym_{\geq 2}(I')$ having no common 2-cycle and such that if $\sigma$ contains a cycle $\omega$ of length $\geq 3$ then $\sigma'$ contains either $\omega$ or $\omega^{-1}$ and vice-versa. 
Denote also by $R^{[2]}_{M_n}(z)$ the sum of the terms  with $\sigma\in\Sym_{\geq 2}(I)$ and $\sigma'\in\Sym_{\geq 2}(I')$ having no common 2-cycle and finally define $R^{[3]}_{M_n}(z):=|P_{M_n{[0]}}(z)|^2$ which is the sum of all terms.
Note that $R^{[i]}_{M_n}(z)$ is a partial sum of $R^{[i+1]}_{M_n}(z)$.

\begin{lemma} \label{lemma_R0}
There is a universal constant $c>0$ such that 
$$|\E(R^{[0]}_{M_n}(z))|\leq cn^{n+1} e^{{1\over 2}z^2-n} \quad \text{for all}\quad  z\in\R.$$
\end{lemma}
\proof
Let $\Sigma(I,I')$ denote the set of $(\sigma,\sigma')$ in $\Sym_2(I)\times \Sym_2(I')$ such that $\sigma$ and $\sigma'$ have no common 2-cycle. From the definition of $R^{[0]}_{M_n}(z)$, we get
$$\E(R^{[0]}_{M_n}(z))=\sum_{(\sigma,\sigma')\in\Sigma(I,I')} (-1)^{|I|+|I'|+\sign(\sigma)+\sign(\sigma')} z^{n-|I|} \overline z^{n-|I'|}.$$
On the other hand, we obtain from the proof of Lemma \ref{lemma_exp_M} that 
$$|\E(P_{M_n[0]}(z))|^2=\E(P_{M_n[0]}(z))\E(P_{M_n[0]}(\overline z))=  \sum_{\sigma\in\Sym_2(I)\atop \sigma'\in\Sym_2(I')} (-1)^{|I|+|I'|+\sign(\sigma)+\sign(\sigma')} z^{n-|I|} \overline z^{n-|I'|}.$$
In order to deduce the lemma from the estimate in Lemma  \ref{lemma_exp_M}, we need to study the contribution of $\sigma$ and $\sigma'$
when they have common 2-cycles.  

Consider a subset $K$ of $\llbracket 1, n \rrbracket$ of length $2k$ and two subsets $J,J'$ of $\llbracket 1,n \rrbracket\setminus K$. Define $I=K\cup J$ and $I'=K\cup J'$. Consider also $\sigma_0\in \Sym_2(K)$, $\sigma_1\in \Sym_2(J)$, $\sigma_1'\in\Sym_2(J')$. Define $\sigma:=\sigma_0\circ\sigma_1$, $\sigma':=\sigma_0\circ \sigma_1'$ and consider the sum
$$T_{K,\sigma_0}:=\sum_{\sigma_1,\sigma_1',J,J'} (-1)^{|I|+|I'|+\sign(\sigma)+\sign(\sigma')} z^{n-|I|} \overline z^{n-|I'|}.$$
This expression represents the sum of terms in $|\E(P_{M_n[0]}(z))|^2$ such that $\sigma$ and $\sigma'$ have $k$ common 2-cycles supported by $K$ which constitute $\sigma_0$. Note that they may have other common cycles supported by $\llbracket 1,n \rrbracket\setminus K$. Using the expansion of $|\E(P_{M_n[K][0]}(z))|^2$, we see that 
$$T_{K,\sigma_0}=|\E(P_{M_n[K][0]}(z))|^2.$$ 

Consider now the sum
$$S_k:=\sum_{|K|=2k}\sum_{\sigma_0,\sigma_1,\sigma_1',J,J'} (-1)^{|I|+|I'|+\sign(\sigma)+\sign(\sigma')} z^{n-|I|} \overline z^{n-|I'|}=\sum_{|K|=2k,\sigma_0} |\E(P_{M_n[K][0]}(z))|^2.$$
This is a combination of terms  in $|\E(P_{M_n[0]}(z))|^2$ such that $\sigma$ and $\sigma'$ have at least  
$k$ common 2-cycles. If for such a term, $\sigma$ and $\sigma'$ have exactly $l$  common 2-cycles, with $l\geq k$, then the term appears in the last double sum 
  ${l\choose k}$ times because this is the number of $(K,\sigma_0)$ for which the considered term appears in $T_{K,\sigma_0}$. 

Note that $S_0=|\E(P_{M_n[0]}(z))|^2$ and $\E(R^{[0]}_{M_n}(z))$ can be written as a linear combination of $S_k$. Namely, we deduce from the above discussion and the identity 
$$\sum_{k=0}^l (-1)^k{l\choose k} =0 \quad \text{for} \quad l\geq 1$$
that
$$\E(R^{[0]}_{M_n}(z))=\sum_{0\leq k\leq n/2} (-1)^k S_k.$$
The number of choices for $K$ of length $2k$ is ${n\choose 2k}$. For each $K$, the number of choices for $\sigma_0$ is $(2k-1)!!$. By Lemma \ref{lemma_exp_M} applied to $M_n[K][0]$ and Stirling's approximation, we have
\begin{eqnarray*}
|S_k| & \leq &  {n\choose 2k} (2k-1)!! c(n-2k)^{n-2k+1} e^{{1\over 2} z^2-n+2k} \\
& = & {cn!(2k-1)!! (n-2k)^{n-2k+1}e^{2k} \over (2k)!(n-2k)! n^{n+1}} n^{n+1} e^{{1\over 2} z^2-n} \\
& = & {cn! (n-2k)^{n-2k+1}e^{2k} \over 2^kk!(n-2k)! n^{n+1}} n^{n+1} e^{{1\over 2} z^2-n} \\
&\leq & {c'\over 2^k k!} n^{n+1} e^{{1\over 2} z^2-n},
\end{eqnarray*}
for some universal constant $c'>0$. The lemma follows easily.
\endproof

\begin{lemma} \label{lemma_l_cycle}
Let $l$ be any integer with $1\leq l\leq n$. 
Then the number of $\sigma\in\Sym(\llbracket 1,n\rrbracket)$ which contain exactly $l$ cycles is less than or equal to $2nn! 2^{-l}$. 
\end{lemma}
\proof
We prove the lemma by induction in $n$. The case $n=1$ is clear. Suppose the lemma holds for $1,\ldots,n-1$. We prove it for $n$. 
The case $l=1$ is clear because the cardinality of $\Sym(\llbracket 1,n\rrbracket)$ is $n!$. We assume that $l\geq 2$.   
Consider $\sigma$ as in the lemma. So 1 belongs to a cycle $\omega$ in $\sigma$. Denote by $m$ the length of this cycle. We have $m\leq n-1$. The number of choices of $\omega$ is the number of cycles of length $m$ containing 1. So it is equal to  $(n-1)\ldots (n-m+1)$ when $m$ is fixed. 

Consider such an  $\omega$ and denote by $K$ its support. 
Observe that $\sigma$ is the composition of $\omega$ with some $\sigma'\in\Sym(\llbracket 1,n\rrbracket\setminus K)$ which has exactly $l-1$ cycles. By induction hypothesis, the number of choices for $\sigma'$ is at most $2(n-m) (n-m)! 2^{-l+1}$ when $\omega$ is fixed. It follows that  the number of $\sigma\in\Sym(\llbracket 1,n\rrbracket)$ which contain exactly $l$ cycles is at most 
$$\sum_{m=1}^{n-1}(n-1)\ldots (n-m+1)2(n-m) (n-m)! 2^{-l+1}=2nn!2^{-l} {2\over n^2} \sum_{m=1}^{n-1} (n-m)\leq 2nn!2^{-l}.$$
The lemma follows.
\endproof

\begin{lemma}
There is a universal constant $c>0$ such that 
$$|\E(R^{[1]}_{M_n}(z))|\leq cn^{n+4} e^{{1\over 2}z^2-n} \quad \text{for all}\quad  z\in\R.$$
\end{lemma}
\proof
Consider $\sigma\in\Sym_{\geq 2}(I)$ and $\sigma'\in\Sym_{\geq 2}(I')$ having no common 2-cycle such that if $\omega$ is a cycle in $\sigma$ of length $\geq 3$ then either $\omega$ or $\omega^{-1}$ is a cycle of $\sigma'$ and vice versa. Let $l$ be the number of these cycles and $K$ the union of their supports. Define $J:=I\setminus K$ and $J':=I'\setminus K$. We can write $\sigma=\sigma_0\circ\sigma_1$ and $\sigma':=\sigma_0'\circ\sigma_1'$ with $\sigma_0,\sigma_0'\in \Sym_{\geq 3}(K)$, $\sigma_1\in\Sym_2(J)$ and $\sigma_1'\in\Sym_2(J')$ such that  if $\omega$ is a cycle of $\sigma_0$  then either $\omega$ or $\omega^{-1}$ is a cycle of $\sigma_0'$ and vice versa.

Using that $\sign(\sigma_0)=\sign(\sigma_0')$, we have
\begin{eqnarray*}
R^{[1]}_{M_n}(z) & = & \sum (-1)^{|J|+|J'|+\sign(\sigma_1)+\sign(\sigma_1')}\xi_{\sigma_0}\xi_{\sigma_0'} \xi_{\sigma_1}\xi_{\sigma_1'} z^{n-|I|} \overline z^{n-|I'|} \\
& = & \sum_{K,\sigma_0,\sigma_0'} \xi_{\sigma_0}\xi_{\sigma_0'}  R^{[0]}_{M_n[K]}.
\end{eqnarray*}
Observe that  $|\E(\xi_{ij}^2)|\leq \E(|\xi_{ij}|^2)=1$ and hence
$|\E( \xi_{\sigma_0}\xi_{\sigma_0'})|\leq 1$. We then deduce from Lemma \ref{lemma_R0}, applied to $R^{[0]}_{M_n[K]}$, that 
$$|\E(R^{[1]}_{M_n}(z))| \leq \sum_{K,\sigma_0,\sigma_0'} c(n-k)^{n-k+1} e^{{1\over 2}z^2-n+k},$$
where we denote by $k$ the length of $K$.

When $K$ and $\sigma_0$ are fixed and $l$ is the number of cycles in $\sigma_0$, there are $2^l$ choices for $\sigma_0'$. When $K$ and $l$ are  fixed, by Lemma \ref{lemma_l_cycle} applied to $k$ instead of $n$, 
the number of choices for $\sigma_0$ is not more than $2kk!2^{-l}$. Finally, the number of choices for $K$ with length $k$ is ${n\choose k}$. We then deduce from the last estimate on $|\E(R^{[1]}_{M_n}(z))|$ that
$$|\E(R^{[1]}_{M_n}(z))| \leq \sum_{0\leq l\leq k\leq n} {n\choose k}2kk!2^{-l}2^lc(n-k)^{n-k+1} e^{{1\over 2}z^2-n+k}.$$
Using Stirling's approximation, we obtain
$$|\E(R^{[1]}_{M_n}(z))| \leq 2c\sum_{0\leq l\leq k\leq n} {n!k(n-k)^{n-k+1}e^{k} \over (n-k)! } e^{{1\over 2}z^2-n} \leq c' n^{n+4} e^{{1\over 2}z^2-n},$$
for some universal constant $c'>0$. The lemma follows.
\endproof

\begin{lemma} \label{lemma_R_2}
There is a universal constant $c>0$ such that 
$$|\E(R^{[2]}_{M_n}(z))|=|\E(R^{[1]}_{M_n}(z))|\leq cn^{n+4} e^{{1\over 2}z^2-n}\quad \text{for all} \quad  z\in\R.$$
\end{lemma}
\proof
The inequality was obtained in the last lemma. So we only need to prove the equality. Observe that 
$$R^{[2]}_{M_n}(z)-R^{[1]}_{M_n}(z)=\sum (-1)^{|I|+|I'|+\sign(\sigma)+\sign(\sigma')}\xi_\sigma \xi_{\sigma'} z^{n-|I|} \overline z^{n-|I'|},$$
where we only consider $\sigma\in\Sym_{\geq 2}(I)$ and $\sigma'\in\Sym_{\geq 2}(I')$ such that $\sigma$ contains a cycle $\omega$ of length $\geq 3$ but $\sigma'$ doesn't contain $\omega$ nor $\omega'$, or vice versa. 

We can represent $\sigma,\sigma'$ as two oriented graphs with vertices $1,2,\ldots,n$. Each of these graphs is composed by disjoint cycles. The above condition on $\sigma,\sigma',\omega,\omega'$ implies that some edge belongs to one graph but not to the other regardless the orientation.  Moreover, this edge belongs to a cycle of length $\geq 3$.   We conclude that $\xi_\sigma \xi_{\sigma'} $ contains a factor $\xi_{ij}$ but not $\xi_{ij}^2$ nor $\xi_{ij}\xi_{ji}$. Thus, the expectation of $\xi_\sigma \xi_{\sigma'} $ vanishes. The lemma follows.
\endproof

\begin{lemma} \label{lemma_R_3}
Let $\beta > 0$ be a number such that $\E(|\xi_{ij}|^4)\leq \beta$ for $i\not=j$. Then
there is a universal constant $c>0$ such that 
$$|\E(R^{[3]}_{M_n}(z))|\leq cn^{n+4} e^{\beta/2}e^{{1\over 2}z^2-n} \quad \text{for all} \quad z\in\R.$$
\end{lemma}
\proof
Recall that 
$$R^{[3]}_{M_n}(z)=\sum_{\sigma\in\Sym_{\geq 2}(I)\atop \sigma'\in\Sym_{\geq 2}(I')} (-1)^{|I|+|I'|+\sign(\sigma)+\sign(\sigma')}\xi_\sigma \xi_{\sigma'} z^{n-|I|} \overline z^{n-|I'|}.$$
Let $K$ denote the union of the supports of common 2-cycles in $\sigma$ and $\sigma'$. Denote by $2k$ its length. So we can write $J:=I\setminus K$, $J':=I'\setminus K$, 
$\sigma=\sigma_0\circ\sigma_1$ and  $\sigma'=\sigma_0\circ\sigma'_1$ for some $\sigma_0\in\Sym_2(K)$, $\sigma_1\in\Sym_{\geq 2}(J)$ and  $\sigma'_1\in\Sym_{\geq 2}(J')$.  We deduce that
$$R^{[3]}_{M_n}(z)=\sum_{K,\sigma_0} \xi_{\sigma_0}^2 R^{[2]}_{M_n[K]}(z).$$

By definition of $\beta$, we have $\E(\xi_{\sigma_0}^2)\leq \beta^k$.
The number of sets $K$ of cardinality $2k$ is ${n\choose 2k}$ and when $K$ is fixed the number of choices for $\sigma_0$ is $(2k-1)!!$. Applying Lemma \ref{lemma_R_2} to $M_n[K]$ instead of $M_n$, together with Stirling's approximation, gives
\begin{eqnarray*}
|\E(R^{[3]}_{M_n}(z))| &\leq & \sum_{k\leq n/2}{n\choose 2k} (2k-1)!! \beta^k c(n-2k)^{n-2k+4} e^{{1\over 2}z^2-n+2k} \\
& \leq & cn^{n+4} e^{{1\over 2}z^2-n}\sum_{k\leq n/2}{n! \beta^k(n-2k)^{n-2k+4}e^{2k}\over 2^k k! (n-2k)! n^{n+4}} \\
& \leq & c'n^{n+4} e^{{1\over 2}z^2-n}\sum_{k\leq n/2}{\beta^k\over 2^k k!} \ \leq \ c'e^{\beta/2} n^{n+4} e^{{1\over 2}z^2-n}, 
\end{eqnarray*}
where $c'>0$ is a universal constant. The lemma follows.
\endproof

The following property is the main estimate in this section.

\begin{proposition} \label{prop_expectation}
Let $\alpha>0$ and  $\beta > 0$ be two numbers such that $\E(|\xi_{ii}|^2)\leq \alpha$  and $\E(|\xi_{ij}|^4)\leq \beta$ for $i\not=j$.
There is a universal constant $c>0$ such that
$$\E(|Q_{W_n}(z)|^2) \leq c n^4 e^{\alpha+\beta/2} e^{n({1\over 2}z^2-1)}$$
for all $z\in\R$. 
\end{proposition}
\proof
It is enough to show that 
$$\E(|P_{M_n}(z)|^2) \leq c n^{n+4} e^{\alpha+\beta/2} e^{{1\over 2}z^2-n}.$$
Lemma \ref{lemma_R_3} says that a similar inequality holds for $M_n[0]$ instead of $M_n$. 

We have the following expansion where the sum is  taken over $\sigma$ and  the subset $I$ of $\llbracket 1,n \rrbracket$
$$P_{M_n}(z)=\sum_{\sigma\in\Sym(I)} (-1)^{|I|+\sign(\sigma)}\xi_\sigma z^{n-|I|}.$$ 
Let $K$ denote the set of fixed points of $\sigma$ and define $J:=I\setminus K$. Then $\sigma$ can be identified with a permutation $\sigma'$ of $J$ and $\xi_\sigma=\xi_K\xi_{\sigma'}$. It follows that
$$P_{M_n}(z)=\sum_{K}  (-1)^{|K|}\xi_K P_{M_n[K][0]}(z).$$
We deduce that 
$$|P_{M_n}(z)|^2=P_{M_n}(z)P_{M_n}(\overline z)=\sum_{K,K'} (-1)^{|K|+|K'|}\xi_K\xi_{K'} P_{M_n[K][0]}(z)P_{M_n[K'][0]}(\overline z).$$

Observe that when $K\not=K'$ the product $\xi_K\xi_{K'}$ contains a factor $\xi_{ii}$ but not $\xi_{ii}^2$. So the expectation of $\xi_K\xi_{K'}$ vanishes in this case. We only need to consider the case where $K=K'$. Observe that $|\E(\xi_K^2)|\leq \alpha^{|K|}$.  
Lemma \ref{lemma_R_3} applied to $M_n[K]$ instead of $M_n$ gives
$$\E(|P_{M_n}(z)|^2) \leq \sum_K \alpha^{|K|}c(n-|K|)^{n-|K|+4} e^{\beta/2} e^{{1\over 2} z^2-n+|K|}.$$ 
Finally, the number of subsets $K$ of  cardinality $k$ is ${n\choose k}$. So we obtain from the last line and Stirling's approximation that
\begin{eqnarray*}
\E(|P_{M_n}(z)|^2) & \leq &  c\sum_{k\leq n} {n\choose k} \alpha^k (n-k)^{n-k+4} e^{\beta/2}e^{{1\over 2} z^2-n+k} \\
& \leq & cn^{n+4} e^{\beta/2}e^{{1\over 2} z^2-n} \sum_{k\leq n} {n! \alpha^k (n-k)^{n-k+4} e^k\over k! (n-k)! n^{n+4}} \\
& \leq & c'n^{n+4} e^{\beta/2} e^{{1\over 2} z^2-n} \sum_{k\leq n} {\alpha^k \over k!} \ \leq \ c'e^{\alpha+\beta/2} n^{n+4} e^{{1\over 2} z^2-n},
\end{eqnarray*}
where $c'>0$ is a universal constant. The proposition follows.
\endproof

\noindent
{\bf End of the proof of Theorem \ref{th_main}.} 
Consider the natural map $\Pi$ from the space of $n\times n$ Hermitian matrices to  the space of polynomials of degree $n$ which associates each matrix with its characteristic polynomial. Consider the natural probability measure on the space of $n\times n$ Hermitian matrices induced by the probabilities given for their entries. Consider also the direct  image of this probability measure by $\Pi$, which is a probability measure on $\Pc_0^n$.  
By Lemma \ref{lemma_pot_sc} and Proposition \ref{prop_expectation}, we only need to apply  Theorems \ref{th_abstract} and \ref{th_abstract_bis} for  $\alpha=2$,  $c_n=cn^4e^{\alpha+\beta/2}$, $K=[-2,2]$ and $l=L=\R$.
\hfill $\square$

\medskip

Note that we can apply the last assertion of Proposition \ref{prop_sup_mass},  Propositions \ref{prop_mass_out_support}  and \ref{prop_endpoint} in order to get better estimates when the interval $I$ is outside the interval $[-2,2]$ or near the endpoints $\pm2$.

\medskip
\noindent
{\bf Proof of Corollary \ref{cor_graph}.} 
Consider the matrix 
$$W_n':=W_n-{p\over \sigma\sqrt{n}} J_n + {p\over\sigma\sqrt{n}}I_n,$$ 
where $J_n$ is the $n\times n$ matrix with all entries equal to 1 and $I_n$ is the identity $n\times n$ matrix. The entries of $M_n':=\sqrt{n} W_n'$ satisfy the properties (1) and (2) in the introduction with 
$$\alpha=0 \quad \text{and}\quad  \beta={p^3+(1-p)^3\over p(1-p)}\leq {2\over p}=o(n).$$
The spectrum of $W_n'$ can be obtained by translate the spectrum of  $W_n-{p\over \sigma\sqrt{n}} J_n$ by ${\sqrt{p}\over \sqrt{1-p}\sqrt{n}}$. Moreover, the number of eigenvalues of the last matrix in any interval differs from the one of the matrix $W_n$ by 0 or $\pm1$, see \cite{TVW}. Therefore, in order to prove the corollary, we can replace $\mu_{W_n}$ by $\mu_{W_n'}$. Now applying Theorem \ref{th_main} to $W_n'$. The probability to have a graph which does not satisfies the estimate in Corollary  \ref{cor_graph} for some interval $I$ is smaller or equal to $An^A e^{1/p} e^{-A^{-1}\delta n}$ for some universal constant $A>0$. The result follows.
\hfill $\square$

\end{document}